\input amstex
\input amsppt.sty
 \magnification=\magstep1
 \hsize=32truecc
 \vsize=22.2truecm
 \baselineskip=16truept
 \NoBlackBoxes
 \nologo
 \TagsOnRight
 
 \def\Z{\Bbb Z}

 \def\l{\left}
 \def\r{\right}
 \def\bg{\bigg}
 \def\({\bg(}
 \def\[{\bg[}
 \def\){\bg)}
 \def\]{\bg]}
 \def\t{\text}
 \def\f{\frac}

 \def\mo{\roman{mod}}

 \def\bi{\binom}
 \def\eq{\equiv}
 
 \def\ls{\leqslant}
 \def\gs{\geqslant}

 \def\bi{\binom}

 \def\Remark{\medskip\noindent{\it  Remark}}
 \def\Ack{\medskip\noindent {\bf Acknowledgment}}
\hbox{Preprint, {\tt arXiv:1312.1166}}
\bigskip
  \topmatter
  \title On $a^n+bn$ modulo $m$\endtitle
  \author Zhi-Wei Sun \endauthor
  \affil Department of Mathematics, Nanjing University\\
 Nanjing 210093, People's Republic of China
  \\  zwsun\@nju.edu.cn
  \\ {\tt http://math.nju.edu.cn/$\sim$zwsun}
\endaffil

 \abstract Let $a$ and $m>0$ be integers. We show that for any integer $b$ relatively prime to $m$, the set  $\{a^n+bn:\ n=1,\ldots,m^2\}$  contains
 a complete system of residues modulo $m$. We also pose several
 conjectures for further research; for example, we conjecture that any integer $n>1$ can be written as 
 $k+m$ with $2^k+m$ prime, where $k$ and $m$ are positive integers.
 \endabstract
 \thanks 2000 {\it Mathematics Subject Classification.}
Primary 11A07, 11B75; Secondary 05A10, 11A41, 11D61, 11P32.
\newline
\indent Supported by the National Natural Science Foundation (grant 11171140) of China.
 \endthanks
 \endtopmatter
  \document
 \heading {1. Introduction}\endheading

 Let $p$ be a prime and let $a$ be a positive integer. In 2011 the author and W. Zhang [SZ] showed that for each $k=p^a,p^a+1,\ldots,2p^a-1$
 the set $\{\bi nk:\, n=0,1,2,\ldots\}$ is dense in the ring of $p$-adic integers,
 i.e., it contains a complete system of residues modulo any power of $p$.

In this paper, we establish the following new result.

 \proclaim{Theorem 1.1} Let $a,\ b$ and $m>0$ be integers. If $b$ is relatively prime to $m$, then
 the set $\{a^n+bn:\ n=1,\ldots,m^2\}$ contains a complete system of residues
 modulo $m$.
 \endproclaim

Our proof of Theorem 1.1 will be given in Section 2.

Now we pose several conjectures for further research.

\proclaim{Conjecture 1.1} For any integers $a$ and $m>0$, the set
$$\{a^n-n:\ n=1,\ldots,2p_m-3\}$$
contains a complete system of residues modulo $m$, where $p_m$ denotes the $m$-th prime. We may also replace $a^n-n$ by $a^n+n$.
\endproclaim
\Remark\ 1.1. For example, $\{2^n-n:\ n=1,\ldots,195\}$ contains a complete system of residues modulo $29$, and $195<2p_{29}-3=2\times 109-3=215$.
\medskip

The following conjecture was motivated by Theorem 1.1 in the cases $b=\pm1$.

\proclaim{Conjecture 1.2} The diophantine equation
   $$x^n+n=y^m\ \t{with}\ m,n,x,y >1$$
has only two integral solutions:
$$5^2+2=3^3\quad\ \t{and}\ \quad 5^3+3=2^7.$$
Also, the diophantine equation
   $$x^n-n=y^m\ \t{with}\ m,n,x,y >1$$
has only two integral solutions:
$$2^5-5 =3^3\quad\ \t{and}\ \quad 2^7-7=11^2.$$
 \endproclaim
\Remark\ 1.2. Conjecture 1.2 seems difficult.
\medskip

By Theorem 1.1, $2^k-k$ or $2^k+k$ modulo a positive integer $m$ behaves better than $2^k$. So our following conjecture is somewhat reasonable.

\proclaim{Conjecture 1.3} {\rm (i)} For any integer $n>1$, there is a positive integer $k<n$ with $n-k+2^k$ prime.
Also, for any integer $n>3$ there is a positive integer $k<n$ with $n+k+2^k$ prime.

{\rm (ii)} Any integer $n>3$ can be written as $p+(2^k-k)+(2^m-m)$, where $p$ is a prime, and $k$ and $m$ are positive integers.
\endproclaim
\Remark\ 1.3. (i) We have verified the first assertion in Conjecture 1.3(i) for $n$ up to $2\times 10^6$ except for
$n=1657977$. For $n=421801$, the least positive integer $k$ with $n-k+2^k$ prime is $149536$.
For $n=1657977$, the least positive integer $k$ with $n-k+2^k$ prime is greater than $2\times10^5$.
We have also verified the second assertion in Conjecture 1.3(i) for all $n=4,5,\ldots,2\times10^6$; for example,
the least positive integer $k$ with $299591+k+2^k$ prime is $51116$.

(ii) The author verified Conjecture 1.3(ii) for all $n=4,5,\ldots,2\times10^8$.
(After learning Conjecture 1.3(ii) from the author, Qing-Hu Hou checked it for all $n$ up to $10^{10}$ without finding any counterexample.)
In contrast, R. Crocker [C] proved in 1971 that
there are infinitely many positive odd integers not of the form $p+2^k+2^m$, where $p$ is a prime, and $k$ and $m$ are positive integers.
See also H. Pan [P] for a further refinement of Crocker's result.
\medskip

Our following conjecture is somewhat similar to Conjectures 1.1 and 1.2.

\proclaim{Conjecture 1.4} {\rm (i)} Let $m$ be any positive integer. Then, either of the sets
$$\{p_n-n:\ n=1,\ldots,2p_m-3\}\quad\t{and}\quad\{np_n:\ n=1,\ldots,2p_m-3\}$$
contains a complete system of residues modulo $m$.

{\rm (ii)} For any non-constant integer-valued polynomial $P(x)$ with positive leading coefficients,
there are infinitely many positive integers $n$ with $p_n-n\ (\t{or}\ p_n+n)$ in the range $P(\Z)$
if and only if $\deg(P)<4$. Also, for any positive integer $n\not=3$, the number $np_n+1$ is not of the form $x^m$
with $m,x\in\{2,3,\ldots\}$.
\endproclaim
\Remark\ 1.4. For example, both $\{p_n-n:\ n=1,\ldots,11\}$ and $\{np_n:\ n=1,\ldots,11\}$ contain a complete system of residues modulo $4$.
Note that $2p_4-3=11$ and that $3p_3+1=3\times5+1=2^4$.

\proclaim{Conjecture 1.5} {\rm (i)} Let $m$ be any positive integer. Then either of the following four sets
$$\align&\l\{\bi{2n}n+n:\ n=1,\ldots,\l\lfloor\f{m^2}2\r\rfloor+3\r\},
\\&\l\{\bi{2n}n-n:\ n=1,\ldots,\l\lfloor\f{m^2}2\r\rfloor+15\r\},
\\&\l\{C_n-n:\ n=1,\ldots,\l\lfloor\f{m^2}2\r\rfloor+7\r\},
\\&\l\{C_n+n:\ n=1,\ldots,\l\lfloor\f{m^2}2\r\rfloor+23\r\}
\endalign$$
contains a complete system of residues modulo $m$, where $C_n$ denotes the Catalan number $\f1{n+1}\bi{2n}n=\bi{2n}n-\bi{2n}{n+1}$.

{\rm (ii)} For any integer $n>2$, neither $\bi{2n}n+n$ nor $\bi{2n}n-n$ has the form $x^m$ with $m,x\in\{2,3,\ldots\}$.
For any integer $n>3$, neither $C_n+n$ nor $C_n-n$ has the form $x^m$ with $m,x\in\{2,3,\ldots\}$.
\endproclaim

\Remark\ 1.5. We also have some other conjectures similar to Conjecture 1.5.

\medskip

For a positive integer $n$, let $p(n)$ denote the number of ways to write $n$ as a sum of positive integers with the order of addends ignored.
Concerning the partition function $p(n)$, M. Newman [N] conjectured that
for any integers $m>0$ and $r$ there are infinitely many positive integers $n$ with $p(n)\eq r\pmod m$.

\proclaim{Conjecture 1.6} {\rm (i)} For any positive integer $n$, we have $p(n)\not= x^m$ for all $m,x\in\{2,3,\ldots\}$.

{\rm (ii)} Any integer $n>3$ can be written in the form $p+p(k)+p(m)$, where $p$ is a prime, and $k$ and $m$ are positive integers.

{\rm (iii)} Each integer $n>4$ can be written as $p+2^k+p(m)$, where $p$ is a prime, and $k$ and $m$ are positive integers.
\endproclaim
\Remark\ 1.6. (i) We also conjecture that no Bell number has the form $x^m$ with $m,x\in\{2,3,\ldots\}$.

(ii) For the representation function corresponding to Conjecture 1.7(ii), see [S, A202650].
\medskip

For a positive integer $n$, let $q(n)$ denote the number of ways to write $n$
as a sum of {\it distinct} positive integers with the order of addends ignored.
The function $q(n)$ is usually called the strict partition function.

\proclaim{Conjecture 1.7} {\rm (i)} For any integer $m>0$ and $r$, there are infinitely many positive integers $n$
with $q(n)\eq r\pmod m$.

{\rm (ii)} For each integer $n>1$, $q(k)q(n-k)+1$ is prime for some $0<k<n$. Also, for any integer $n>5$, $q(k)q(n-k)-1$ is prime for some
$0<k<n$.

{\rm (iii)} For any integer $n>1$, there is a positive integer $k<n$ such that $p(k)^2+q(n-k)^2\ ($ or $p(k)+q(n-k))$ is prime.
\endproclaim
\Remark\ 1.7. It is known that
$$p(n)\sim\f{e^{\pi\sqrt{2n/3}}}{4\sqrt3 n}\ \ \t{and}\ \ q(n)\sim\f{e^{\pi\sqrt{n/3}}}{4(3n^3)^{1/4}}\ \quad\t{as}\ n\to+\infty$$
(cf. [HR] and [AS, p.\,826]).
Part (i) of Conjecture 1.7 is an analogue of Newman's conjecture on the partition function, for example, the least positive integer $n$
with $q(n)\eq31\ (\mo\ 42)$ is $8400$.
Part (ii) implies that there are infinitely many primes $p$ with $p-1$ a product of two strict partition numbers.
Similarly, part (iii) implies that there are infinitely many primes of the form $p(k)^2+q(m)^2$ with $k$ and $m$ positive integers.
We have verified part (iii) for $n$ up to $10^5$. For some sequences related to parts (ii)-(iii),
see A233417, A232504, A233307, A233346 of [S].

\proclaim{Conjecture 1.8} {\rm (i)} Any integer $n>9$ can be written as $k+m$ with $k$ and $m$ positive integers such that
$2^{\varphi(k)/2+\varphi(m)/6}+3$ is prime. Also, any integer $n>13$ can be written as $k+m$ with $k$ and $m$ positive integers such that
$2^{\varphi(k)/2+\varphi(m)/6}-3$ is prime.

{\rm (ii)} Any integer $n>25$ can be written as $k+m$ with $k$ and $m$ positive integers such that
$3\times 2^{\varphi(k)/2+\varphi(m)/8}+1$ is prime. Also, any integer $n>14$ can be written as $k+m$ with $k$ and $m$ positive integers such that
$3\times 2^{\varphi(k)/2+\varphi(m)/12}-1$ is prime.
\endproclaim
\Remark\ 1.8. We have verified this for $n$ up to $50000$. The conjecture
implies that there are infinitely many primes in any of the four forms $2^n+3$, $2^n-3$, $2^n3+1$ and $2^n3-1$.

\heading{2. Proof of Theorem 1.1}\endheading

\noindent{\it Proof of Theorem 1.1}. (i) We first use induction to show the claim that
if $m$ is relatively prime to $ab$ then $\{a^n+b n:\ n=1,\ldots,m^2\}$ contains a complete system of residues modulo $m$.

The claim holds trivially for $m=1$.

Now let $m>1$ be relatively prime to $ab$, and assume the claim for smaller values of $m$
relatively prime to $ab$. Write $m=p_1^{a_1}\ldots p_r^{a_r}$, where
$p_1<\ldots<p_r$ are distinct primes and $a_1,\ldots,a_r$ are
positive integers. Note that $m_0=m/p_r$ is relatively prime to $ab$.

Let $r$ be any integer. By the induction hypothesis, there is a
positive integer $k\ls m_0^2$ such that $a^k+b k=r+m_0q_0$ for some $q_0\in\Z$.
Since $b\prod_{i=1}^r(p_i-1)$ is relatively prime to $p_r$, there is
a nonnegative integer $q<p_r$ such that
$$q\times b\prod_{i=1}^r(p_i-1)\eq -q_0\pmod{p_r}.$$
Set $n=k+m_0q\prod_{i=1}^r(p_i-1)$. Since $\varphi(m)$ divides $m_0\prod_{i=1}^r(p_i-1)$, by applying Euler's theorem we obtain
$$a^n+b n\eq a^k +b\(k+m_0q\prod_{i=1}^r(p_i-1)\) =r +m_0\(q_0+bq\prod_{i=1}^r(p_i-1\)\eq r\pmod m.$$
Note that
$$\align 0<k\ls n\ls &m_0^2+m_0(p_r-1)\prod_{i=1}^r(p_i-1)
\\<&m_0^2+m_0(p_r-1)m=m_0^2(1+p_r^2-p_r)<m_0^2p_r^2=m^2.
\endalign$$
This concludes the induction step.

(ii) Now we assume that $b$ is relatively prime to $m$.
If $a$ is also relatively prime to $m$, then the desired result follows from the claim in (i).
Below we suppose that $a$ is not relatively prime to $m$. Write $m=uv$, where $u>1$ and $v>0$ are integers such that
$a$ is divisible by any prime divisor of $u$ and $a$ is relatively prime to $v$.
Let $r$ be an arbitrary integer. As $b$ is relatively to $u$, $bs\eq r\ (\mo\ u)$ for some $s\in\{0,1,\ldots,u-1\}$.
Choose $a_*\in\Z$ with $aa_*\eq1\ (\mo\ v)$. As $a^u$ and $a_*^sbu$ are both relatively prime to $v$, by (i)
there is a positive integer $k\ls v^2$ such that
$$(a^u)^k+a_*^sbuk\eq a_*^s(r-bs)\pmod v.\tag2.1$$
Set $n=uk+s$. Then
$$n\ls uv^2+u-1<u(v^2+1)\ls u(v^2+v^2)\ls u^2v^2=m^2.$$
In view of (2.1),
$$a^{uk+s}+buk\eq r-bs\pmod v,\ \t{i.e.},\ a^n+bn\eq r\pmod v.\tag2.2$$
For any prime divisor $p$ of $u$, the $p$-adic order of $u$ is smaller than $u$ since $p^u\gs 2^u\gs u+1$.
Therefore,
$$a^n+bn=a^{uk+s}+b(uk+s)\eq 0 + bs \eq r\pmod u.\tag2.3$$
Combining (2.2) and (2.3) we obtain that $a^n+bn\eq r\ (\mo\ m)$.

In view of the above, we have completed the proof of Theorem 1.1. \qed

\Ack. The author would like to thank Dr. H. Pan for helpful comments.


\widestnumber \key{SZ}
\Refs

\ref\key AS \by M. Abramowitz and I. A. Stegun (eds.), \book Handbook
of Mathematical Functions with Formulas, Graphs, and Mathematical
Tables, \publ 9th printing, New York, Dover, 1972\endref

\ref\key C\by R. Crocker\paper On a sum of a prime and two powers of two\jour Pacific J. Math.\vol 36\yr 1971\pages 103--107\endref

\ref\key HR\by G. H. Hardy and S. Ramanujan\paper
 Asymptotic formulae in combinatorial analysis, \jour Proc. London Math. Soc. \vol 17\yr 1918\pages 75--115\endref

\ref\key N\by M. Newman\paper Periodicity modulo $m$ and divisibility properties of the partition function
\jour Trans. Amer. Math. Soc.\vol 97\yr 1960\pages 225--236\endref

\ref\key P\by H. Pan\paper On the integers not of the form $p+2^a+2^b$\jour Acta Arith.\vol 148\yr 2011\pages 55--61\endref

\ref\key S\by Z.-W. Sun\paper {\rm Sequences A202650, A231201, A231725, A232398, A232504, A232616, A232862,
A233307, A233346 and A233417 in OEIS (On-Line Encyclopedia of Integer Sequences)}
\jour {\tt http://oeis.org}\endref

\ref\key SZ\by Z.-W. Sun and W. Zhang\paper Binomial coefficients and the ring of $p$-adic integers\jour Proc. Amer. Math. Soc.
\vol 139\yr 2011\pages 1569--1577\endref

\endRefs
\enddocument